# Extending Bricard Octahedra

Gerald D. Nelson
nelso229@comcast.net

**Abstract**

We demonstrate the construction of several families of flexible polyhedra by extending Bricard octahedra to form larger composite flexible polyhedra. These flexible polyhedra are of genus 0 and 1, have dihedral angles that are non-constant under flexion, exhibit self-intersections and are of indefinite size, the smallest of which is a decahedron with seven vertexes.

## 1. Introduction

Flexible polyhedra can change in spatial shape while their edge lengths and face angles remain invariant. The first examples of such polyhedra were octahedra discovered by Bricard [1] in 1897. These polyhedra, commonly known as Bricard octahedra, are of three types, have triangular faces and six vertexes and have self-intersecting faces. Over the past century they have provided the basis for numerous investigations and many papers based in total or part on the subject have been published.

An early paper published in 1912 by Bennett [2] investigated the kinematics of these octahedra and showed that a prismatic flexible polyhedra (polyhedra that have parallel edges and are quadrilateral in cross section) could be derived from Bricard octahedra of the first type. Lebesgue lectured on the subject in 1938 [3]. The relationship between flexible prismatic polyhedra and Bricard octahedra was described in more detail in 1943 by Goldberg [4]. The well known counter-example to the polyhedra rigidity conjecture [5] was constructed by Connelly in 1977 using elements of Bricard octahedra to provide flexibility. A 1990 study [6] by Bushmelev and Sabitov described the configuration space of octahedra in general and of Bricard octahedra specifically. They cite from [3] where it was shown that "every position of one cap (of a Bricard octahedron) completely determines the position of the other half of the octahedron". Two recent papers published in 2002 and 2009, by Stachel and Alexandrov respectively; provide an alternate proof of the flexibility of Bricard octahedra of the third type and prove that the Dehn invariant is a constant during flexure of any Bricard octahedra [7, 8].

Because of the potential application to the construction of various physical mechanisms the subject has also been studied in a more utilitarian context. Papers by Baker [9] and [10] published in 1980 and 1995 are good examples. The first paper has numerous physical examples based upon Bricard octahedra of the first and second types. The latter considers the motion of the third type of Bricard octahedra.

In this paper we describe several families of flexible polyhedra that can be constructed by extending Bricard octahedra to form larger composite flexible polyhedra of both genus 0 and genus 1. The polyhedra in these families are of indefinite size, the smallest being a decahedron with seven vertexes and a hendecahedron with eight vertexes. As a basis for the construction approach, Bricard octahedra are separated into five sub-types based upon the geometric characteristics of the unique caps that can be found in the octahedra.

A polyhedron in this paper is considered to be a closed surface in $\mathbb{R}^3$ formed by simply-connected polygons (the faces) that are joined together at polygon edges of equal length with two faces joined at an edge. Edges are terminated by vertexes that are shared by three or more edges and their associated faces. The sub-set of faces that share a common vertex form an open connected surface that can be completely traversed by moving from face to face by crossing edges and without touching the common vertex. All polyhedra described in this paper are non-convex and have some faces that intersect with one another. Open polyhedra, in which some faces are missing, occur only as intermediate entities in the discussion of the construction of larger closed polyhedra.

Flexible polyhedra are continuously flexible and have faces that are unchanged while their dihedral and solid angles vary over some range of values. For each flexible polyhedron it is possible to identify a geometric parameter, typically a dihedral angle, as the variable of flexion that varies continuously over some range of values.

The term construction is used to describe either the combination of analytical and logical arguments that describe the formation of a polyhedron in mathematical terms or to describe the actual creation, generally computer based, of specific polyhedra. With regard to both meanings, construction is focused on methods to compute the position of vertexes and is of a recursive nature in which vertexes associated with a specific stage of the construction are defined with reference to those of the previous stage.

Finally we note that the flexibility of polyhedra described in this paper is ultimately derived from Bricard's work and this fact is used to limit the scope of this paper. In a more technical sense then the scope is limited to constructions whose flexibility derives from polyhedra with vertexes of index = 4. Additionally, we do not consider the construction of closed flexible polyhedra that do not exhibit self-intersections, thus all closed flexible polyhedra described herein have zero oriented volume; and flexible prismatic polyhedra are not considered.

## 2. Bricard Octahedra – types and sub-types

Three distinct types of flexible octahedra were developed by Bricard based upon an analysis of the geometric properties of the vertexes of an octahedron. These properties describe the angular relationships that exist at tetrahedral vertexes and the behavior when such vertexes are deformed. In this section we summarize relevant details of Bricard octahedra and develop an alternative classification of these octahedra that will be used in subsequent sections to extend these octahedra into larger flexible polyhedra.

As a basis for the following discussion we will use the nomenclature established in the Model Octahedron shown in Figure 1. The vertex labeling conventions shown there will be used throughout the remainder of this paper.

In this model (Fig. 3) the outsides of faces are shown in dark shades; insides in light. Self-intersections show up as transitions between light and dark shades whereas edges and vertexes are shown with finite dimension cylinders and spheres. Faces are designated by a clockwise convention thus the view of the face **$X_1B_0C_0$** is of the inside of the face while the partial view of **$X_0C_0B_0$** is of the outside. The direction of the vector cross product **$(B_0-X_0) \times (C_0-X_0)$** defines the outward pointing normal of the face.

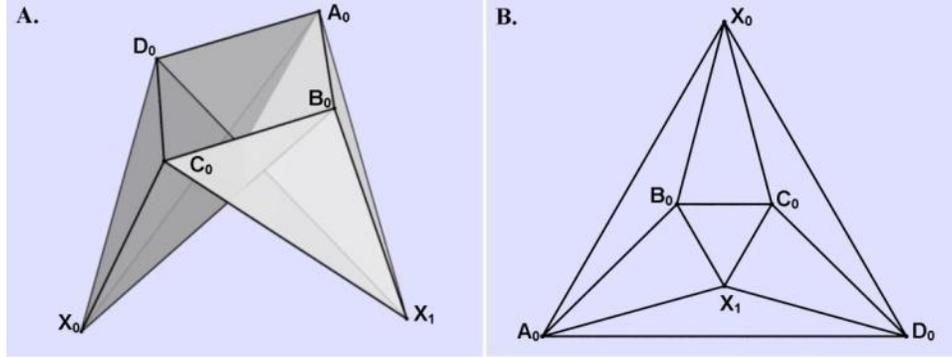

**Figure 1.** Octahedron Model. **(A)** Octahedron Model. **(B)** Octahedron Diagram.

Vertex labeling in this model is not the same as used in [1] and has been adopted here since it generalizes easily into an indexed notation that is used in the description of the construction of polyhedra later in this paper. (Polyhedra images in this figure and all other figures have been created with ray trace software [12].)

Bricard [1] demonstrated the existence of three types of flexible octahedra, the first two of which exhibit axial and planar symmetric motion while the third has two positions in which all vertexes lie in a plane. Alternately these three types can be described by geometric characteristics (Table I) that are more useful in the construction approach that is used in this paper.

|     | Characteristic |
| --- | --- |
| OEE | Opposite Edges Equal |
| AEE | Adjacent Edges Equal |
| OAE | Opposite Angles Equal |
| OAS | Opposite Angles Supplementary |

**Table I.** Geometric Characteristics

The first two characteristics describe edge length relationships that are applicable to Bricard octahedra of the first two types while the angular relationships refer to cap face angles and are applicable to the third type. In the first two types of Bricard octahedra there are six pairs of equal length edges. The edges in the first type are arranged so that all equal length pairs are on opposite edges (edges **$D_0A_0$** and **$B_0C_0$** for example), thus are all OEE. In the second, only four pairs of equal length edges are located opposite of one another, are OEE, while the other two equal length pairs are connected at common vertexes and form a quadrilateral on which there are two pairs of adjacent edges with equal lengths, thus are considered AEE. In the third type, two vertexes, opposite of one another, have opposite face angles that are supplementary, are OAS, while the remaining four vertexes have opposite face angles that are equal, thus are OAE. This configuration of face angles has been previously described [10, Fig. 5].

These characteristics lead naturally to a classification into five unique sub-types that are based upon a mapping of geometric characteristics of the octahedra onto a specific four faced polyhedral surface, or cap, of the octahedron. Generally the specification of a cap is sufficient to define the complete octahedron and provides a convenient method of identification in later sections of this paper. These sub-types are designated as I-OEE, II-AEE, II-OEE, III-OAE and III-OAS. The Roman numeral refers to the type of the associated Bricard octahedron; the abbreviated identifiers, as defined in Table I, refer to specific geometric characteristics of the cap. The cap of interest is the cap at vertex $\mathbf{X_0}$ as illustrated by the Octahedron Model shown in Figure 1. Additionally for all Bricard types, the opposite vertex, $\mathbf{X_1}$, is of the same sub-type as $\mathbf{X_0}$. The defining geometric characteristics of each of these sub-types are summarized in Table II.

| Sub-type | Defining Geometric Characteristics |
|---|---|
| I-OEE | Edge lengths at the cap base are OEE. |
| II-AEE | Edge lengths at the cap base are AEE. |
| II-OEE | Edge lengths at the cap base are OEE; edge lengths at the cap vertex are AEE. |
| III-OAE | Face angles at the cap vertex are OAE. |
| III-OAS | Face angles at the cap vertex are OAS. |

**Table II.** Sub-type Characteristics

While Table II provides a unique definition for each sub-type several other relationships are needed to complete the picture in general and when the cap is considered to be part of an octahedron. Additionally there are octahedra that can satisfy more than one definition; specifically there can be ambiguity between sub-types I-OEE, II-AEE and II-OEE and sub-types III-OAE and III-OAS. For these cases we used the relevant III designation.

The edge length relationships for the first three types are shown in Table III. For the sub-type II-AEE the choice of which of the edges on the base of the cap ($\mathbf{A_0B_0}$, $\mathbf{B_0C_0}$, $\mathbf{C_0D_0}$, $\mathbf{D_0A_0}$) are designated as equal length is by convention.

| Sub-type | Edge Length Relationships |
|---|---|
| I-OEE | $|\mathbf{X_1A_0}| = |\mathbf{X_0C_0}|$, $|\mathbf{X_1B_0}| = |\mathbf{X_0D_0}|$, $|\mathbf{X_1C_0}| = |\mathbf{X_0A_0}|$, $|\mathbf{X_1D_0}| = |\mathbf{X_0B_0}|$, $|\mathbf{C_0D_0}| = |\mathbf{A_0B_0}|$, $|\mathbf{D_0A_0}| = |\mathbf{B_0C_0}|$. |
| II-AEE | $|\mathbf{X_1A_0}| = |\mathbf{X_0C_0}|$, $|\mathbf{X_1B_0}| = |\mathbf{X_0B_0}|$, $|\mathbf{X_1C_0}| = |\mathbf{X_0A_0}|$, $|\mathbf{X_1D_0}| = |\mathbf{X_0D_0}|$, $|\mathbf{B_0C_0}| = |\mathbf{A_0B_0}|$, $|\mathbf{D_0A_0}| = |\mathbf{C_0D_0}|$. |
| II-OEE | $|\mathbf{X_0C_0}| = |\mathbf{X_0A_0}|$, $|\mathbf{X_0D_0}| = |\mathbf{X_0B_0}|$, $|\mathbf{C_0D_0}| = |\mathbf{A_0B_0}|$, $|\mathbf{D_0A_0}| = |\mathbf{B_0C_0}|$, $|\mathbf{X_1C_0}| = |\mathbf{X_1A_0}|$, $|\mathbf{X_1D_0}| = |\mathbf{X_1B_0}|$. |

**Table III.** Edge Length Relationships for sub-types I-OEE, II-AEE II-OEE.

Face angle relationships for the last two sub-types are shown in Table IV. For the sub-type III-OAE the choice of which of the vertexes on the base of the cap ($A_0$, $B_0$. $C_0$, $D_0$) are designated as OAS is by convention.

| Sub-type | Vertex Face Angle Relationships |
|---|---|
| III-OAE | OAE - $X_0$, $X_1$, $A_0$, $C_0$; OAS - $B_0$, $D_0$. |
| III-OAS | OAS - $X_0$, $X_1$; OAE - $A_0$, $B_0$, $C_0$, $D_0$. |

**Table IV.** Face Angle Relationships for sub-types III-OAE and III-OAS.

The five-way classification of the "configuration spaces" of polyhedral caps described by Bushmelev and Sabitov [6] is not directly related to the five sub-types identified here. In that paper the deformations of a cap are classified parametrically and then mapped to Bricard octahedra types as appropriate. For example, Bricard octahedra of the third type are all of the fifth type of "configuration space". In this paper, Bricard octahedra of the third type map into two separate sub-types.

## 3. Recursive Construction

In this section we describe a recursive approach for constructing larger flexible polyhedra by extending Bricard octahedra that is applicable to any of the cap sub-types that are defined in Table II. To illustrate this approach we restrict any diagrams that represent flexible octahedra to those in which the four faces of the primary polyhedral cap, the cap at vertex $X_0$, (Fig. 1) are free of self-intersections. This is done to simplify the diagrams and is not a limitation of the approach. Assuming that the octahedron in Figure 1 is flexible, we introduce the notion of extending flexible octahedra into larger flexible polyhedra by observing that if the edges at the vertex $X_1$ are extended beyond the vertexes of the polyhedron (Fig. 2.A), four quadrilateral surfaces can be created. Edges of the cap at vertex $X_1$ are extended to vertexes $A_1$, $B_1$, $C_1$ and $D_1$ by using the same scale factor on each of the edges.

The four quadrilateral faces, taken together with the triangular faces of the cap at vertex $X_0$, form a composite open flexible octahedron (Fig. 2B) that may or may not have self-intersections. The open flexible octahedral surface is defined by faces {$X_0B_0A_0$, $X_0C_0B_0$, $X_0D_0C_0$, $X_0A_0D_0$, $A_0B_0B_1A_1$, $B_0C_0C_1B_1$, $C_0D_0D_1C_1$, and $D_0A_0A_1D_1$} and results from dropping the four faces {$X_1D_0A_0$, $X_1A_0B_0$, $X_1B_0C_0$, $X_1C_0D_0$} at the vertex $X_1$ from consideration. Also, the direction associated with the added quadrilateral faces (dark faces) is inverted from the direction of the triangular faces at the vertex $X_1$ (light colored) that have been discarded.

Since the lengths of the edge extensions do not affect the flexibility of this composite open octahedron, this method of creating a flexible open octahedron can be repeated to produce successively larger open flexible polyhedra. The extension to form an open composite flexible dodecahedron can be accomplished by creating (Figs. 3A and 3B) a larger flexible octahedron and then extending one of its caps. Edge length extensions from vertex $X_1$ (Fig. 2) have been constructed by scaling so that octahedron consisting of the polyhedral caps at vertexes $X_1$(extended) and $X_2$ is flexible. The octahedral surface in Fig. 3A is augmented by faces {$A_1B_1B_2A_2$, $B_1C_1C_2B_2$, $C_1D_1D_2C_2$, $D_1A_1A_2D_2$} that are formed from the extension of the cap at $X_2$.

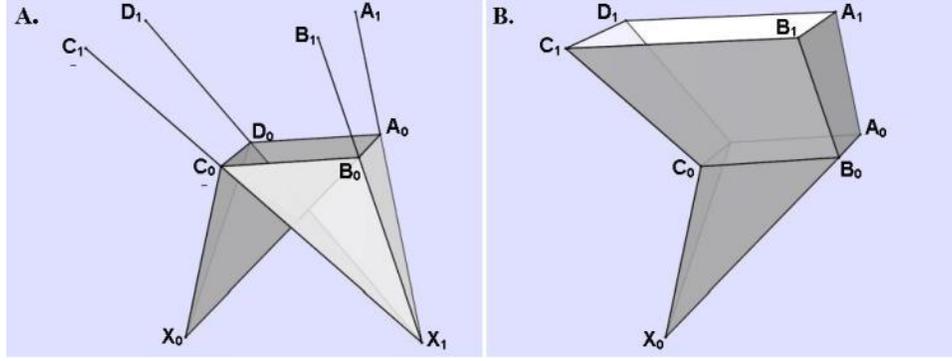

**Figure 2.** Extension of a Flexible Octahedra. **(A)** Extended Vertexes. **(B)** Composite Flexible Octahedral Surface.

Alternately, to produce a closed decahedron exhibiting self-intersections, the quadrilateral extensions can be ignored and the faces of the cap at vertex $X_2$ used to complete the construction (Fig. 3C). The polyhedron is closed by ignoring the extended vertexes $A_2$, $B_2$, $C_2$ and $D_2$ and forming triangular faces with the cap at $X_2$.

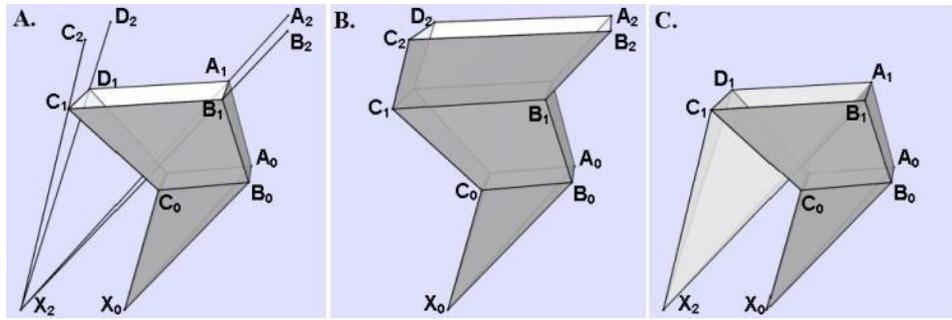

**Figure 3.** Extension of a Flexible Octahedral Surface. **(A)** Extended Vertexes. **(B)** Composite Open Dodecahedron. **(C)** Composite Closed Dodecahedron.

This approach can be formalized into a recursive construction by a consideration of an indefinite number of flexible octahedra that are constructed in stages; the octahedron at a particular stage is constructed from the octahedron of the previous stage and which have a shared flexible surface of quadrilateral faces as described in Figure 4. Two stages of construction are depicted; the i-th, the extension of octahedron $P_{i-1}$ to $P_i$ and the (i+1)st, the extension from of octahedron $P_i$ to $P_{i+1}$ where the octahedra are defined as follows:

$$P_{i-1} = \{A_{i-1}X_{i-1}B_{i-1},\ B_{i-1}X_{i-1}C_{i-1},\ C_{i-1}X_{i-1}D_{i-1},\ D_{i-1}X_{i-1}A_{i-1},$$
$$B_{i-1}X_iA_{i-1},\ C_{i-1}X_iB_{i-1},\ D_{i-1}X_iC_{i-1},\ A_{i-1}X_iD_{i-1}\},$$
$$P_i = \{A_iX_iB_i,\ B_iX_iC_i,\ C_iX_iD_i,\ D_iX_iA_i,$$
$$B_iX_{i+1}A_i,\ C_iX_{i+1}B_i,\ D_iX_{i+1}C_i,\ A_iX_{i+1}D_i\} \text{ and}$$
$$P_{i+1} = \{A_{i+1}X_{i+1}B_{i+1},\ B_{i+1}X_{i+1}C_{i+1},\ C_{i+1}X_{i+1}D_{i+1},\ D_{i+1}X_{i+1}A_{i+1},$$
$$B_{i+1}X_{i+2}A_{i+1},\ C_{i+1}X_{i+2}B_{i+1},\ D_{i+1}X_{i+2}C_{i+1},\ A_{i+1}X_{i+2}D_{i+1}\}.$$

At the i-th stage, assuming that $P_{i-1}$ is flexible, the lengths of the extensions, $|A_i-A_{i-1}|$, $|B_i-B_{i-1}|$, $|C_i-C_{i-1}|$ and $|D_i-D_{i-1}|$ are determined by scaling with the same multiplicative factor so that $P_i$ is flexible. Similarly the extensions to the cap at $X_{i+1}$ can be made so that $P_{i+1}$ is flexible. Thus the surface defined by the quadrilateral faces

and
$$\{A_{i-1}B_{i-1}B_iA_i,\ B_{i-1}C_{i-1}C_iB_i,\ C_{i-1}D_{i-1}D_iC_i,\ D_{i-1}A_{i-1}A_iD_i\}$$
$$\{A_iB_iB_{i+1}A_{i+1},\ B_iC_iC_{i+1}B_{i+1},\ C_iD_iD_{i+1}C_{i+1},\ D_iA_iA_{i+1}D_{i+1}\}$$

that remain when the vertex caps at $X_i$ and $X_{i+1}$ are discarded is also flexible since they are parts of flexible polyhedra.

These surfaces provide a flexible extension to the four triangular faces of the cap at vertex $X_0$ and at the i-th stage create an open flexible polyhedron that may or may not exhibit self-intersections. If at stage i+1 of the construction, the four triangular faces of the cap at vertex $X_{i+2}$ are retained, rather than continuing the construction, a closed flexible polyhedron with self-intersecting faces is created.

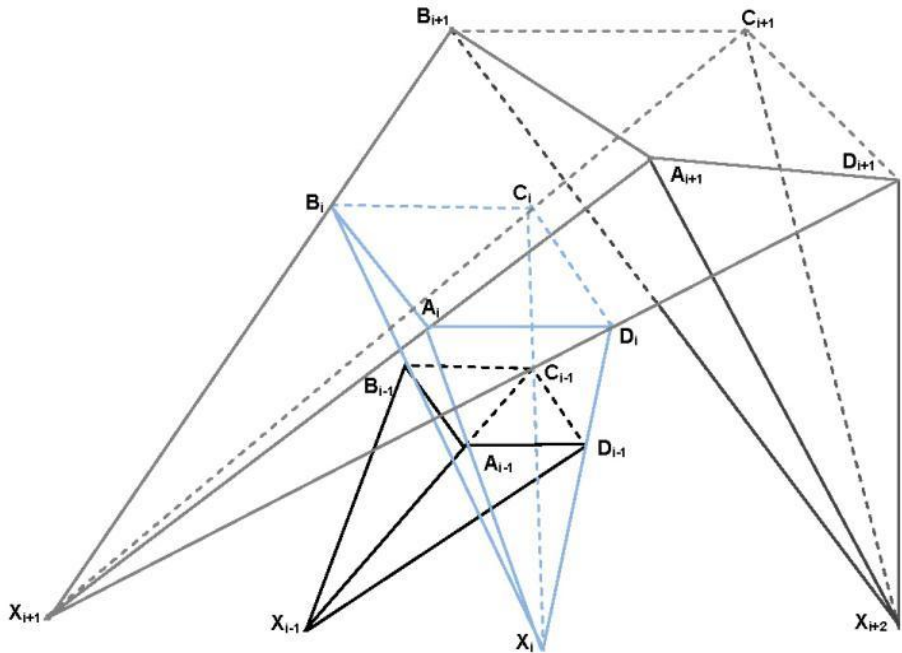

**Figure 4.** Flexible Polyhedra Construction Stages.

The preceding discussion constitutes an inductive proof of the following theorem:

**Theorem 1:** There exist flexible closed polyhedra of genus 0 that have dihedral angles that are non-constant under flexion, exhibit self-intersections and have 4n+2 vertexes (all of index=4) and 4(n+1) faces, eight that are triangular and 4(n-1) that are quadrilateral, for all n>1.

From inspection of Figure 4, it is apparent that the sense of the construction extension can be reversed. Rather than extending one cap of an octahedron to obtain a flexible surface consisting of four quadrilateral faces, quadrilateral faces can be constructed directly by taking some fraction of the cap. The construction process illustrated in Figure 4 is one in which the size of the polyhedron increases; when the process is reversed, the size decreases. Additionally the process can be reversed at any stage.

Finally, it is noted that flexibility of the resultant closed polyhedron can only be achieved by use of Bricard octahedra in the intermediate stages of the construction even though it is true that at any stage, the edges of the intermediate flexible open polyhedron can be arbitrarily extended to yield an open flexible polyhedron. However, the resultant polyhedron cannot be extended, or closed, and remain flexible since this would imply the existence of flexible octahedra different from those of Bricard.

## 4. Genus 1 Construction

We can generalize the recursive method of construction so that the resulting flexible polyhedron is a torus. To illustrate, we consider the construction of the specific flexible polyhedron shown in Figure 5. The small caps at vertexes $X_1$ and $X_2$, are dropped from consideration during the first and second stages of construction. The cap at vertex $X_3$ is retained during the third and last stage.

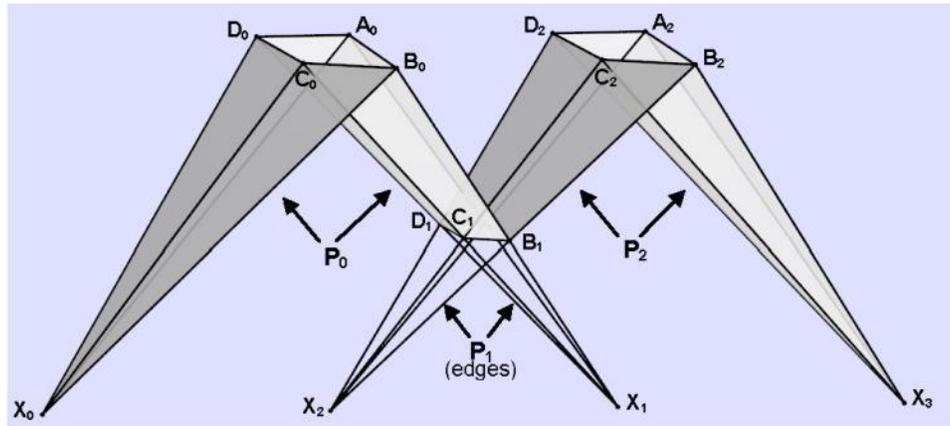

**Figure 5.** Hexadecahedron and construction extensions.

This hexadecahedron is constructed from three identical flexible octahedra; two of which are the same size and the third of which is of reduced size by a scale factor. The intermediate flexible octahedra used in this construction are:

$P_0 = \{A_0X_0B_0, B_0X_0C_0, C_0X_0D_0, D_0X_0A_0, B_0X_1A_0, C_0X_1B_0, D_0X_1C_0, A_0X_1D_0\}$,
$P_1 = \{B_1X_1A_1, C_1X_1B_1, D_1X_1C_1, A_1X_1D_1, A_1X_2B_1, B_1X_2C_1, C_1X_2D_1, D_1X_2A_1\}$ and
$P_2 = \{A_2X_2B_2, B_2X_2C_2, C_2X_2D_2, D_2X_2A_2, B_2X_3A_2, C_2X_3B_2, D_2X_3C_2, A_2X_3D_2\}$.

Edge lengths of the octahedron $P_1$ are computed by scaling the edges of $P_0$ by a factor $f<1$. Edge lengths of the octahedron $P_2$ are computed by scaling the edges of $P_1$ by $1/f$, resulting in two identical octahedra $P_0$ and $P_2$.

From inspection it is apparent that a fourth flexible octahedron can be constructed by extending the edges at the two caps $X_0$ and $X_3$ of the hexadecahedron to intersect at vertexes $A_3$, $B_3$, $C_3$ and $D_3$ (Fig. 6A). From the symmetry of the construction the necessary scale factor, F say, is seen to be $F = |X_3-X_0|/|X_1-X_0|$. But $|X_3-X_0|=2|X_1-X_0|-|X_2-X_1|$ and $f = |X_2-X_1|/|X_1-X_0|$ so that $F = 2-f$.

Removal of the original caps at $X_0$ and $X_3$ and addition of the eight faces associated with the vertexes $A_3$, $B_3$, $C_3$ and $D_3$, results in the sixteen facetted torus shown in Figure 6B. The choice of direction associated with the faces that have been added to create the torus is based on the direction associated with faces that are retained from the underlying hexadecahedron. This choice is a construction option since the direction of all faces could be completely reversed and be equally valid.

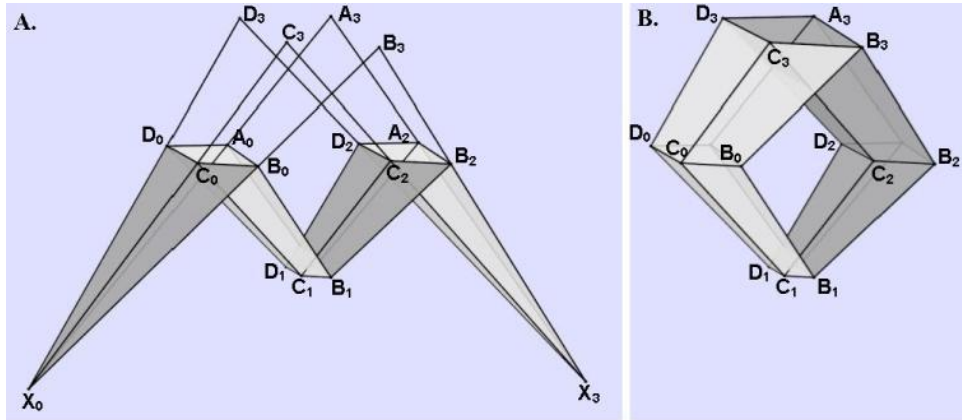

**Figure 6.** Torus Construction. **(A)** Hexadecahedron and extensions. **(B)** Torus with sixteen faces.

This construction approach can be generalized to produce tori of virtually any size and shape. By taking M identical Bricard octahedra arranged as shown in Figure 5 and overlapping so that there are N levels of smaller octahedra an intermediate structure can be created that presents innumerable construction opportunities. An example structure that illustrates this for M=5 and N=4 is shown in Figure 7A. Since the structure is comprised of flexible octahedra the positions of vertexes at the intersections of the various octahedra are unchanged under flexion. By retaining the quadrilaterals that appear between connected intersections, such as those that are labeled in Figure 7A, a variety of tori can be constructed. An example is shown in Figure 7B.

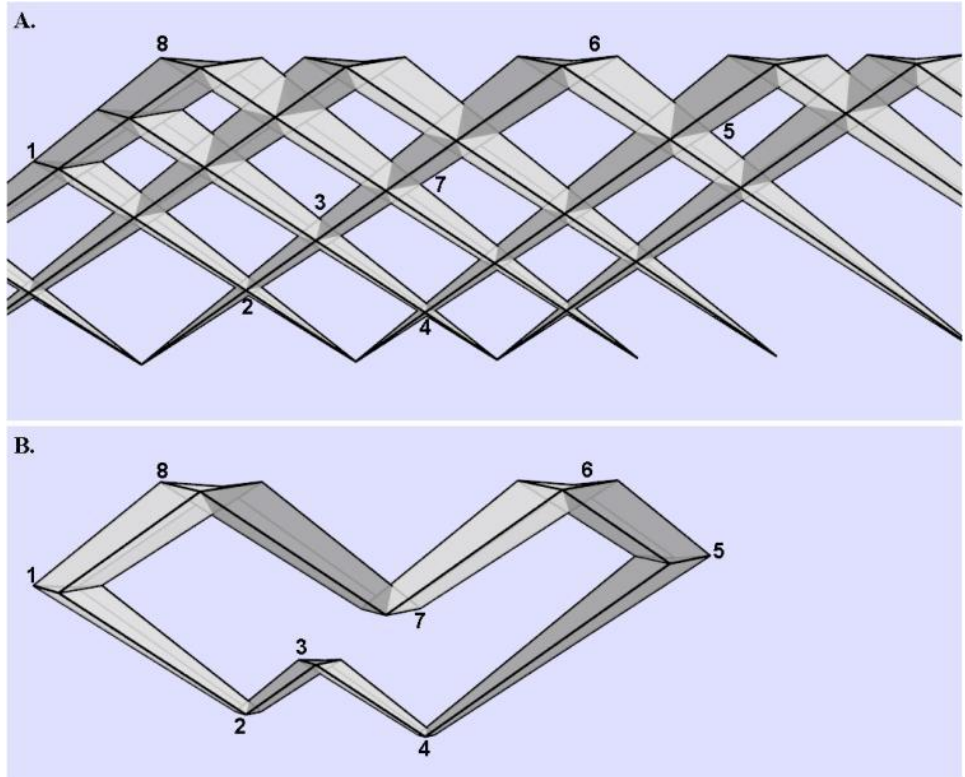

**Figure 7.** General Torus Construction. **(A)** Multiple Overlain Flexible Octahedra. **(B)** Torus with thirty-two faces.

Since M and N are indefinite in size we can state the following theorem:

**Theorem 2:** There exist flexible closed polyhedra of genus 1 that have dihedral angles that are non-constant under flexion, exhibit self-intersections and have 8n vertexes, all of index=4, and an equal number of quadrilateral faces, for n>1.

## 5. Some Generalizations

Five distinct families of flexible polyhedra that are associated with the sub-type identifiers I-OEE, II-AEE, II-OEE, III-OAE and III-OAS can be constructed by applying the recursive methods described above to Bricard octahedra. They consist of polyhedra that satisfy Theorems 1 and 2 and are constructed entirely from intermediate octahedra of the same sub-type. Genus 0 polyhedra have two tetrahedral caps at opposite ends of the polyhedra that are separated from one another by quadrilateral faces that are arranged in annular bands of four faces per band. Each of these bands represents the contribution of one construction stage of the composite polyhedron. Genus 1 polyhedra are comprised solely of quadrilateral faces which are arranged in annular bands of four faces per band.

Additionally there are a number of variations of the constructions of genus 0 polyhedra that can be considered as contributing a sixth family of polyhedra. At each stage of the construction three variations can occur; the intermediate polyhedra that are used in successive stages can be created by methods other than scaling, there can be a change of sub-type in the transition from one stage to the next or fewer than four edge length extensions can be used in a stage. The resulting polyhedra are characterized by enumeration of the specific sub-types and variations that are used in the successive stages of construction and are open ended in that any number of such variations can be used.

Extending the intermediate polyhedra by methods other than scaling is accomplished on a stage by stage basis and consists of the use of straight-forward, sub-type specific computations in which the magnitude of some of the edge length extensions are treated as parameters and other lengths computed accordingly given that the face angles at the cap vertex are fixed. One application of this approach is to specify the length of two extensions at each stage of construction and then compute the remaining two extension lengths.

Intermixing of sub-types can be done between I-OEE and II-AEE sub-types and between II-OEE and III-OAE sub-types. In both of these cases the cap face angle configurations are compatible in the sense that it possible to construct caps that lead to either sub-type when the face angles at the cap vertex are known. Other sub-type combinations have unique cap face angle configurations that are not compatible. At the i-th stage of this sort of construction the sub-type that is associated with octahedron $P_{i-1}$ (Fig. 4) is different than the sub-type that is associated with $P_i$ and the two caps defined by vertexes {$X_i$, $A_{i-1}$, $B_{i-1}$, $C_{i-1}$, $D_{i-1}$} and the vertexes {$X_i$, $A_i$, $B_i$, $C_i$, $D_i$} are of different sub-types.

The use of fewer than four edge length extensions during stages of the construction permits the addition of some triangular faces rather than quadrilateral faces and in some cases the addition of fewer that four faces. This does not apply to II-OEE constructions as two non-zero parameters are always required with this sub-type. For sub-types I-OEE and II-AEE, one parameter can be set to zero; thus construction with three non-zero edge length extensions is possible. For sub-types III-OAE and III-OAS, up to two parameters can be set to zero; thus construction with either two or three non-zero edge length extensions is possible. Also, in some cases, III-OAE constructions can be performed with three parameters set to zero provided that the assignment of OAS vertexes is handled appropriately.

While a full exposition of these possibilities is beyond the scope of this paper we will illustrate the latter case with sufficient constructions to prove the following theorem:

**Theorem 3:** There exist flexible closed polyhedra of genus 0 that have dihedral angles that are non-constant under flexion, exhibit self-intersections, have two positions in which all vertexes are co-planar and have n faces, for all n>9.

Specifically we show explicit constructions of a decahedron with seven vertexes and a hendecahedron with eight vertexes based upon sub-type III-OAE. Thus, since we have previously established a mechanism for extending flexible polyhedra to any size by systematically adding four quadrilaterals at a stage, we argue that the theorem is proved.

In addition to providing support for the proof of Theorem 3, the decahedron provides an example of a flexible polyhedron having the smallest number of faces that can be constructed using the methods described in this paper. This construction (illustrated in Fig. 8) is accomplished by making the extension from the vertex $X_1$ (Fig. 2A) with three of the extension lengths effectively set to zero. Also the parameterization of the polyhedron $P_0$ is such that the extended vertex, $A_1$ in this case, is on the surface of the cap at $X_1$ (Fig. 8C).

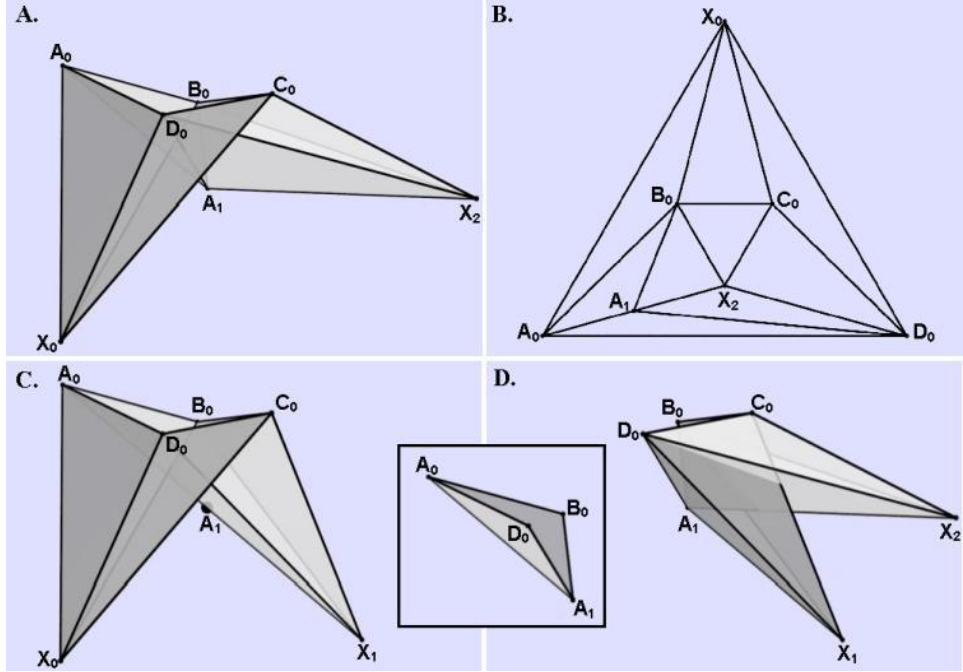

**Figure 8.** III-OAE Flexible Decahedron. **(A)** Decahedron. **(B)** Decahedron Diagram. **(C)** Octahedron $P_0$. **(D)** Octahedron $P_1$. **Inset.** Extension of cap at $X_1$.

As part of the construction of the decahedron (Fig. 8A) two intermediate octahedra, $P_0$ (Fig. 8C) and $P_1$ (Fig. 8D), are constructed and the decahedron formed from a composite of the cap at $X_0$ in $P_0$, the extension of the cap at $X_1$ of $P_1$(the two triangular faces shown in the inset, $A_0A_1D_0$ and $A_0B_0A_1$) and the cap at $X_2$ in $P_1$.

Construction of these flexible octahedra is accomplished by a parameterization of the octahedra followed by edge length and face angle computations. Parameters for $P_0$ are taken to be five independent parameters of the two adjacent faces $X_0B_0A_0$ and $X_0C_0B_0$. (Construction of flexible octahedra, of the third Bricard type, is known [2, Sec. 29] to entail five independent parameters.) Additionally, $P_0$ is constructed with vertexes $B_0$ and $D_0$ treated as OAS. $P_1$ is parameterized by two of the faces, $X_1C_0B_0$ and $X_1D_0C_0$, on the cap at $X_1$, based upon completion of the construction of $P_0$ and the vertexes $A_1$ and $C_0$ are treated as OAS. As a consequence, the construction of $P_1$ generates a new cap at $X_1$, a cap that has the vertex $A_0$ extended to $A_1$.

Edge length and face angle computations that complete the definition of a parameterized octahedron are based upon application of the law of sines and of application of equations derived by Bricard [1, Eqs. 4 and 5]. These define the behavior of OAE and OAS vertexes under deformation. We observe that the parameterization of the two adjacent faces of an III-OAE cap at vertex $X_0$ leaves two angles on the other two faces of the cap, $X_0A_0D_0$ and $X_0D_0C_0$, unresolved; the angles $\angle X_0A_0D_0$ and $\angle D_0C_0X_0$ for example. Further, since $A_0$ and $C_0$ are OAE vertexes $\angle X_1A_0B_0 = \angle X_0A_0D_0$ and $\angle B_0C_0X_1 = \angle D_0C_0X_0$. The following discussion defines equations that relate these two angles to one another.

By applying the law of sines to the faces $X_1A_0B_0$ and $X_1B_0C_0$ and eliminating the edge length $|X_1B_0|$ from the resulting equations it is easy to see that the two angles of interest are related by:

$$\operatorname{ctn}\Gamma_2 = a \operatorname{ctn}B_1 + b, \tag{1}$$

where

$$a = L_2 \cos\Gamma_1 / L_1 \sin B_2 \text{ and}$$

$$b = (L_2 \sin\Gamma_1 - L_1 \cos B_2) / L_1 \sin B_2,$$

with edge lengths and angles defined by $L_1 = |A_0B_0|$, $L_2 = |B_0C_0|$, $B_1 = \angle X_1A_0B_0$, $B_2 = \angle X_1B_0C_0$, $\Gamma_1 = \angle A_0B_0X_1$ and $\Gamma_2 = \angle B_0C_0X_1$.

Applying equations derived by Bricard [1, Eqs. 4 and 5] to the vertexes $A_0$, $B_0$ and $C_0$ one finds that

$$\operatorname{ctn}(\Gamma_2/2) = k \operatorname{ctn}(B_1/2), \tag{2}$$

where

$$k = \tan(\beta_1/2) / \tan(\gamma_2/2).$$

with angles $\beta_1$ and $\gamma_2$ defined by $\beta_1 = \angle B_0A_0X_0$ and $\gamma_2 = \angle X_0C_0B_0$.

Combining Eqs. 1 and 2 yields the quadratic equation:

$$(k - a)\operatorname{ctn}^2(\Gamma_2/2) - 2bk \operatorname{ctn}(\Gamma_2/2) + k(ak - 1) = 0. \tag{3}$$

With $\Gamma_2$ known from the solution of Eq. 3 and $B_1$ from Eq. 1 the remaining faces angles are known directly from the OAS and OAE assignments for the vertexes and the entire octahedron is fully defined. The construction of $P_1$ is completed in the same manner with Eq. 1 applied to faces $X_2B_0C_0$ and $X_2C_0D_0$ and Eq. 2 applied to vertexes $B_0$, $C_0$ and $D_0$.

It is noted that there other descriptions of the construction of Bricard octahedra of the third type that are based upon a more geometric approach and which may be of interest; see [1, Sec. 11; 2, Sec.28; 7, Sec. 1; 8, Sec.4]. We use the above approach as it seems to be computationally simpler and easier to apply with the recursive constructions described in this paper.

Hendecahedron construction (Fig. 9) is a simple variation on the above decahedron construction. The $P_0$ polyhedron (Fig. 9C) is constructed as described above while for $P_1$ (Fig. 9D) the edge extension associated with the edge $X_1B_0$ is set to some to non-zero value as part of the construction sequence. This results in the creation of three faces (Fig. 9 Inset) during the first stage of construction, two triangles, $B_1B_0C_0$ and $D_0A_0A_1$, and one quadrilateral, $A_0B_0B_1A_1$.

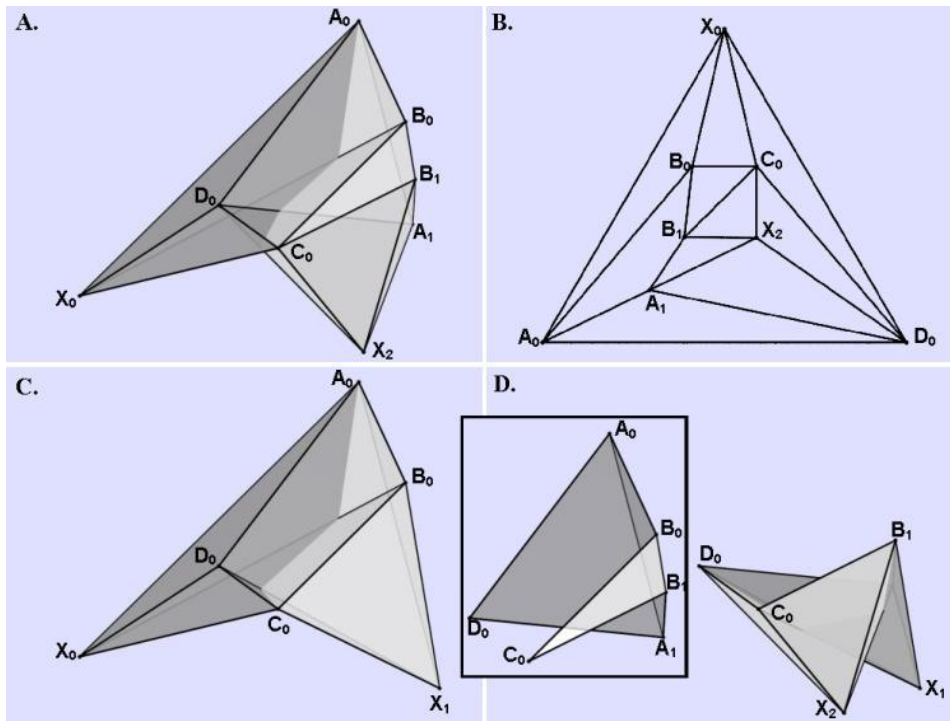

**Figure 9.** III-OAE Flexible Hendecahedron. **(A)** Hendecahedron. **(B)** Hendecahedron Diagram. **(C)** Octahedron $P_0$. **(D)** Octahedron $P_1$. **Inset.** Extension of cap at $X_1$.

## 6. Construction Details

During the preparation of this paper we have created a number of flexible polyhedra (computer constructions) that support the production of figures, to test construction algorithms and simply to observe an interesting and beautiful four dimensional phenomenon. Here we provide some comments about computational techniques that have been found to be useful. Additionally, Appendix A provides a description, parameters and graphic images, of the polyhedra that are used in previous figures and of some other example constructions.

As we have seen, our construction approach for both genus 0 and genus 1 flexible polyhedra depends upon the ability to construct a series of related intermediate Bricard octahedra and to create composite polyhedra from selected portions of the intermediate constructions. Generally the treatment of initial polyhedron **P$_0$** (Fig. 4) in this series is somewhat different than the treatment of the octahedra that occur in subsequent stages. In principle, the construction of all intermediate octahedra is carried out to completion; however the means of characterization is somewhat different. Construction of polyhedron **P$_0$** entails the choice of an appropriate variable of flexion and coordinate model, and sufficient parameterization to support the computation of the vertex positions; whereas subsequent octahedra may be characterized by a simple scale factor.

By convention, we treat the dihedral angle along the edge **X$_0$A$_0$** (Fig. 1) as the variable of flexion, thus all position coordinates are effectively a function of this variable.

Characterization of I-OEE, II-AEE and II-OEE sub-types depends upon six parameters which are nominally taken to be edge lengths. For the first two sub-types, the edge lengths of the cap at **X$_0$** are sufficient to completely parameterize the octahedron, for II-OEE sub-types two additional lengths from the cap at **X$_1$** are required. The characterization of all edge lengths for the complete octahedron is shown in Table III.

Characterization of III-OAE and III-OAS sub-types depends upon five parameters that we nominally take to be some combination of the face angles and edge lengths of the faces **X$_0$B$_0$A$_0$** and **X$_0$A$_0$D$_0$**. In this case, all face angles at the vertex **X$_0$** are known from the parameterization and it is readily apparent that one additional parameter, say the edge length |**X$_0$C$_0$**| or the face angle ∠**C$_0$B$_0$X$_0$** is required to complete the construction of the cap at **X$_0$** and of the octahedron. This angle can be computed as described in the derivation of Eq. 3. However, all parameterzations of the faces **X$_0$B$_0$A$_0$** and **X$_0$A$_0$D$_0$** do not result in valid values for the face angle ∠**C$_0$B$_0$X$_0$**, thus some experimentation and analysis is frequently necessary to find suitable parameters. For the sub-type, III-OAE two possible solutions may exist, depending upon which pair of vertexes **A$_0$** and **C$_0$** or **B$_0$** and **D$_0$** is designated to be OAS.

For purposes of display the sub-type specific coordinate models shown in Table V are convenient.

| I-OEE | **X$_0$**(0,0,z$_z$), **X$_1$**(0,0,-z$_z$), **A$_0$**(a$_x$,a$_y$,a$_z$) , **B$_0$**(b$_x$,b$_y$,b$_z$), **C$_0$**(-a$_x$,a$_y$,-a$_z$), and **D$_0$**(-b$_x$,b$_y$,-b$_z$). |
|---|---|
| II-AEE | **X$_0$**(0,z$_y$,z$_z$), **X$_1$**(0,z$_y$,-z$_z$), **A$_0$**(0,a$_y$,a$_z$), **B$_0$**(b$_x$,b$_y$,0), **C$_0$**(0,a$_y$,-a$_z$), and **D$_0$**(d$_x$,d$_y$,0). |
| II-OEE | **X$_0$**(z$_x$,z$_y$,0), **X$_1$**(x$_x$,x$_y$,0), **A$_0$**(0,a$_y$,a$_z$) , **B$_0$**(0,b$_y$,b$_z$), **C$_0$**(0,a$_y$,-a$_z$), and **D$_0$**(0,b$_y$,-b$_z$). |
| III-OAE, III-OAS | **X$_0$**(0,0,z$_z$), **X$_1$**(x$_x$,x$_y$,x$_z$), **A$_0$**(0,0,a$_z$), **B$_0$**(b$_x$,0,b$_z$), **C$_0$**(c$_x$,c$_y$,c$_z$), and **D$_0$**(d$_x$,d$_y$,0). |

**Table V**. Sub-type Coordinate Models.

I-OEE coordinate model emphasizes the axial symmetry (y-axis) of octahedra of this sub-type. II-AEE and II-OEE coordinate models emphasize the planar symmetry (x-y plane) of these octahedra and in both cases the four co-planar vertexes of the octahedron are located in the y-z plane and move symmetrically about the origin. III-OAE and III-OAS coordinate models fix the face **X$_0$B$_0$A$_0$** in the x-z plane and align the variable of flexion with the z-axis; thus motion of the vertexes $C_0$, $D_0$, and $X_1$ is circular and is defined by rotations about the edges $X_0B_0$, $X_0A_0$ and $A_0B_0$ respectively. This coordinate model is also suitable for the I-OEE, II-AEE and I-OEE sub-types however the motion symmetry associated with these sub-types is lost. In all cases these coordinates can be readily computed from parameters and the variable of flexion, by rotation, triangulation or directly from length equations. The II-AEE and II-OEE coordinate models support the specific edge length relations defined in Table III.

The construction of intermediate octahedra at subsequent stages requires that one cap of the previous stage be extended, or contracted, in such a manner that it is part of another flexible octahedron. This can be done by scaling all edge lengths of the cap being extended by the same scale factor and resulting in a series of octahedra that are successive enlargements or contractions of the $P_0$ octahedron. Scaling is applicable to constructions for all sub-types and provides the choice of a single parameter at each construction stage.

Two edge length parameters can be independently specified at each stage for the sub-types I-OEE and II-AEE. The remaining two edge lengths can then be obtained from the solution of the quadratic equations that result from applying the law of cosines to edge lengths at vertex $X_i$. For both sub-types, two edges that are opposite of one another can be treated as independent parameters and for sub-type I-AEE, two edges that are adjacent to one another at vertex $X_i$ can be treated as independent parameters.

Three edge length parameters can be independently specified for the sub-types III-OAE and III-OAS. The construction is identical to the construction described above for the initial octahedron and is applied to the faces $X_iB_iA_i$ and $X_iA_iD_i$ which are completely determined at the i-th stage of construction. The undetermined angle at the vertex $B_i$, the angle $\angle C_iB_iX_i$, is determined from the solution of the relevant equations depending upon the face angle configuration at the vertex $B_i$.

## 7. Conclusion

In this paper we have 1) developed a classification of Bricard octahedra into five sub-types that are based upon the geometric characteristics of one of the caps of an octahedron, 2) developed a recursive method of construction that makes it possible to extend Bricard octahedra into larger flexible polyhedra including flexible tori, 3) identified six families of extended composite flexible polyhedra, and 4) created example constructions for each of the six families. While Bricard octahedra have been widely studied and remarked upon it appears that these results are novel.

## Acknowledgements

The author wishes to thank Dr. Kari M. Nelson and Ms. Heidi Nelson Ries for their editorial comments on earlier drafts of this paper and his wife Verla Nelson for her encouragement and support.

**Gerald Nelson** is a retired software engineer; worked for Honeywell, Inc. (when it was still a Minnesota based company) and MTS Systems Corp. of Eden Prairie, Minnesota. He has a Masters degree in Mathematics from the University of Minnesota. His interest in flexible polyhedra stems from reading an article [11] about the Bellows Conjecture.

# Appendix A

This appendix describes polyhedra used as illustrations in the body of this paper along with several other example polyhedra constructions. Each polyhedron is described by the parameters used for construction and by some images of the spatial shape and motion of the polyhedron under flexion. In all cases the flexibility of the polyhedron being described has been numerically tested by computation of several invariant quantities: oriented volume, Dehn invariant, edge lengths, face areas and the sum of the two solid angles at the initial and closing caps of the polyhedron. The latter sum is $= 4\pi$ for all sub-types and in the case of polyhedra of the third type each solid angle is individually $= 2\pi$. Error or variations in these quantities over the full range of flexion are always observed to be insignificant when using 15-digit floating point precision. Typically, when dealing with edge lengths ~10, errors that are of the order of $10^{-13}$ are observed.

Although parameter values that are used in the construction of the following examples have been rounded to whole or half numbers these are nominal values in the sense that there is no particular significance associated with a specific given value. Generally any given parameter value can be varied continuously over a reasonable interval with no significant change in the nature of the construction. In the following examples, length parameters are unit-less.

The Model Octahedron (Fig. 1A) and the initial stages of the recursive construction (Figs. 2 and 3) are illustrated using a sub-type I-OEE octahedron. Here we continue the extensions (Fig. 3B) to construct a flexible polyhedron with 32 faces (Fig. 10).

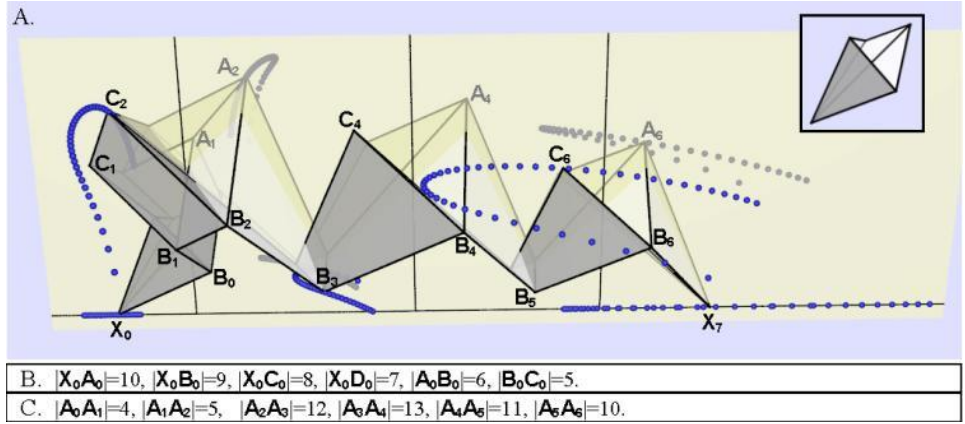

**Figure 10.** I-OEE Flexible Polyhedron. **A.** 32 Facetted Polyhedron. **B.** Cap Construction Parameters. **C.** Extension Parameters. **Inset.** Cap at $X_0$.

The axial motion that is characteristic of Bricard octahedra of the first type is carried into this construction (Fig. 10) as illustrated by the motion traces (blue spheres). These show that the two vertexes, $X_0$ and $X_7$, move in a straight line while the remaining vertexes move symmetrically in a pair wise fashion about axes that are orthogonal to the straight line. The motion trace of $X_7$ appears to be irregular since the motion is reversed in direction at some point during flexion. The light shaded plane contains the line of motion and three of the axes of symmetry. The pair of vertexes $C_6$ (in front of the plane) and $A_6$ (behind the plane) move symmetrically about the third axis. Parameters (Fig. 10C) for this example are edge length extensions that at each stage of construction are converted into scale factors that are applied to all four edges to achieve the specified length. The construction strategy that is employed relies on extensions that both increase and decrease the size of polyhedron. The first two stages, the fourth and the sixth, extend the intermediate caps beyond the octahedron while the third and fifth extend onto the intermediate caps. The seventh stage closes the polyhedron.

An example of an II-AEE flexible polyhedron construction is shown in Figure 11.

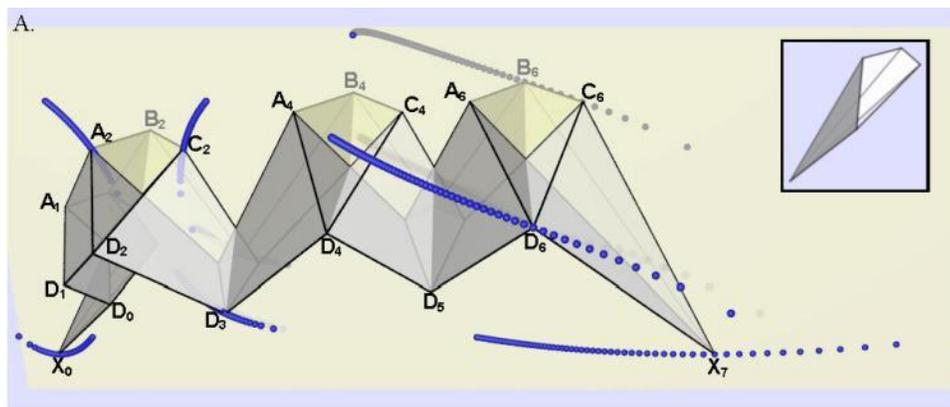

**Figure 11.** II-AEE Flexible Polyhedron. **A.** 32 Facetted Polyhedron. (cont.)

| B. | $|X_0A_0|=10$, $|X_0B_0|=11$, $|X_0C_0|=12$, $|X_0D_0|=13$, $|A_0B_0|=3$, $|C_0D_0|=8$. |
| C. | $|A_0A_1|=4$, $|A_1A_2|=5$, $|A_2A_3|=12$, $|A_3A_4|=13$, $|A_4A_5|=11$, $|A_5A_6|=10$. |

**Figure 11.** II-AEE Flexible Polyhedron. **B.** Cap Construction Parameters. **C.** Extension Parameters. **Inset.** Cap at $X_0$.

The same extension parameters and extension strategy as used in the I-OEE example (Fig. 10C) are used with different cap parameters (Fig. 11B). The planar motion of Bricard octahedra of the second type appears in this construction. The vertexes $X_0$, $X_1$, $A_i$ and $C_i$ for i=1..6 are co-planar and move in the light shaded plane. The pairs of vertexes $B_i$ and $D_i$ move symmetrically with respect to this plane. The vertexes $B_i$ are behind the plane, $D_i$ are in front.

Planar motion of Bricard octahedra of the second type takes on a somewhat different form in the II-OEE Flexible Icosahedron shown in Figure 12. In it, the vertexes $X_0$ and $X_4$ move in a plane (the light shaded plane), while each of sets of four vertexes $A_i$, $B_i$, $C_i$ and $D_i$ for i=1..4 are co-planar and move symmetrically in planes that are orthogonal to the plane of motion of $X_0$ and $X_4$.

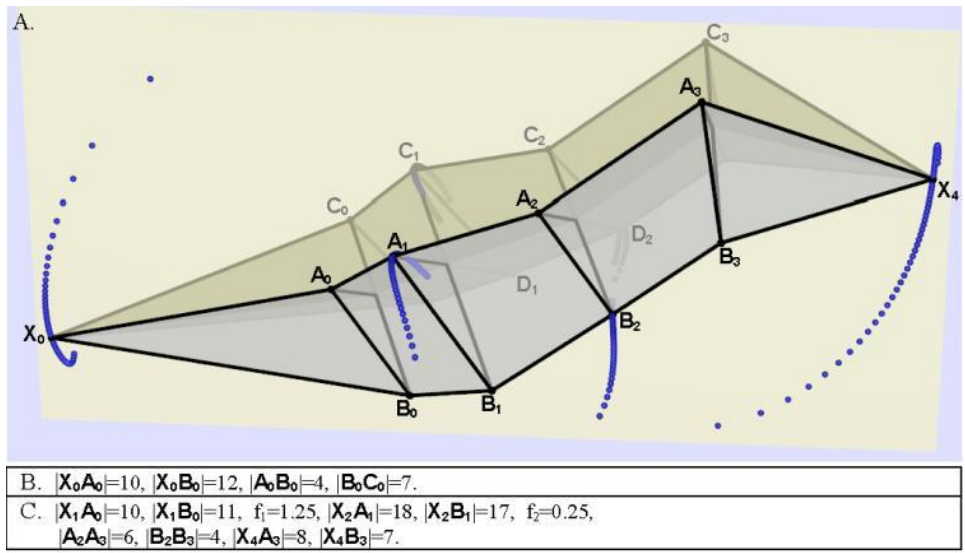

| B. | $|X_0A_0|=10$, $|X_0B_0|=12$, $|A_0B_0|=4$, $|B_0C_0|=7$. |
| C. | $|X_1A_0|=10$, $|X_1B_0|=11$, $f_1=1.25$, $|X_2A_1|=18$, $|X_2B_1|=17$, $f_2=0.25$, $|A_2A_3|=6$, $|B_2B_3|=4$, $|X_4A_3|=8$, $|X_4B_3|=7$. |

**Figure 12.** II-OEE Flexible Polyhedron. **A.** Icosahedron. **B.** Cap Construction Parameters. **C.** Extension Parameters.

The construction strategy (Fig. 12) at each stage consists of a specification of the edge lengths of intermediate caps along with a factor of the cap edge length to retain for the stage. Since the four vertexes $A_i$, $B_i$, $C_i$ and $D_i$ are co-planar the intermediate cap edge lengths are essentially un-constrained by the previous stage. Additionally, extensions that have parallel edges ($A_2A_3$, $B_2B_3$, $C_2C_3$, $D_2D_3$ for example) are possible.

III-OAE construction is illustrated by the flexible hecatohedron shown in Figure 13. This polyhedron is the result of a 24 stage construction in which there are three edge length deltas specified at each stage. A new type three Bricard octahedron is created at each stage. The intermediate edges are extended so that successively larger faces are created resulting in the large cap at the vertex **X₂₄** where the construction is terminated.

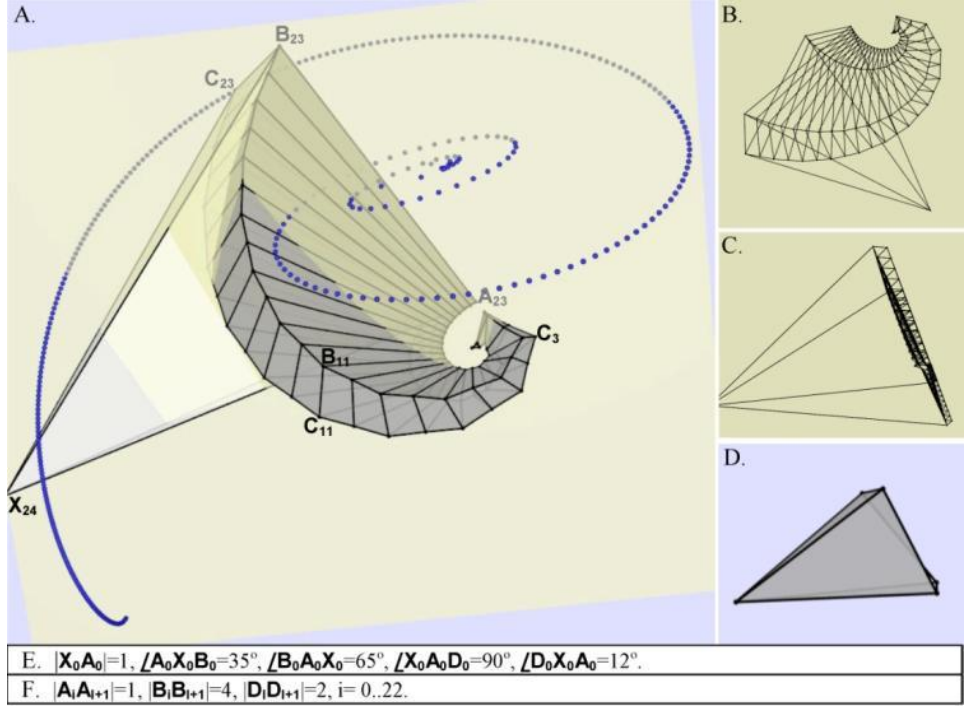

E. $|X_0A_0|=1$, $\angle A_0X_0B_0=35°$, $\angle B_0A_0X_0=65°$, $\angle X_0A_0D_0=90°$, $\angle D_0X_0A_0=12°$.
F. $|A_iA_{i+1}|=1$, $|B_iB_{i+1}|=4$, $|D_iD_{i+1}|=2$, i= 0..22.

**Figure 13.** III-OAE Flexible Polyhedron. **A.** Hecatohedron. **B.** Structure – first flat position. **C.** Second flat position. **D.** Cap at **X₀**. **E. P₀** Cap Construction Parameters. **F.** Extension Parameters.

The intermediate octahedra used in this hecatohedron are created with the same construction strategy that is represented by Eqs. 1, 2 and 3 with a different choice when the Bricard equations [1, Eqs. 4 and 5] are applied to the vertexes **A_i**, **B_i** and **C_i**. These choices lead to a variation in the form of Eqs. 2 and 3:

$$\mathrm{ctn}(\Gamma_2/2) = k \tan(B_1/2),  \quad\quad\quad (2a)$$

where

$$k = \mathrm{ctn}(\beta_1/2) / \tan(\gamma_2/2).$$

with angles $\beta_1$ and $\gamma_2$ defined by $\beta_1 = \angle B_iA_iX_i$, $\gamma_2 = \angle X_iC_iB_i$, $B_1 = \angle X_{i+1}A_iB_i$ and $\Gamma_2 = \angle B_iC_iX_{i+1}$.

Combining Eqs. 1 and 2a yields the quadratic equation:

$$(k + a)\operatorname{ctn}^2(\Gamma_2/2) - 2bk\operatorname{ctn}(\Gamma_2/2) - k(1 + ak) = 0. \qquad (3a)$$

The motion of the vertexes of this polyhedron is quite dramatic as indicated by the motion trace shown for vertex $C_{23}$ that shows the motion as the vertex moves from one first flat position (Fig. 13B) to the second flat position (Fig. 13D). The hecatohedron (Fig. 13A) is displayed at a position where the dihedral angle at edge $X_0A_0$ is $= 90°$.

III-OAS construction is illustrated by the flexible hecatohedron shown in Figure 14.

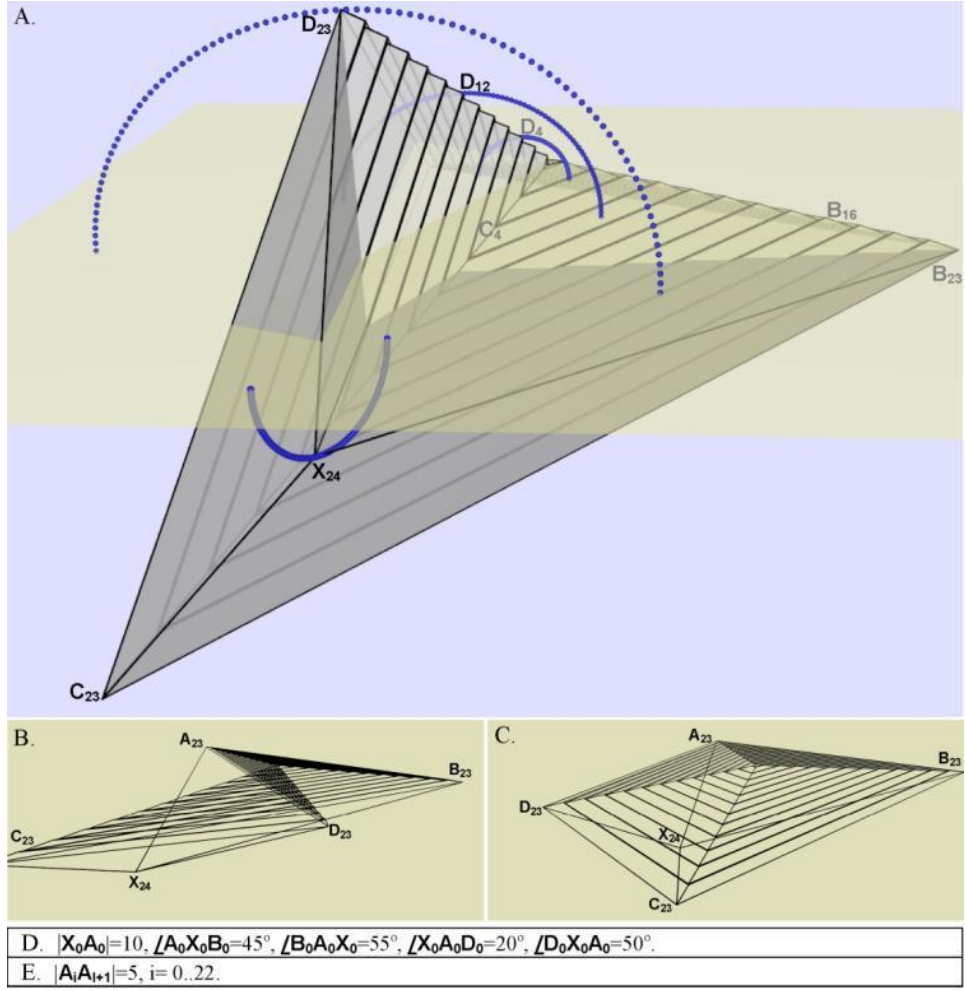

**Figure 14.** III-OAS Flexible Polyhedron. **A.** Hecatohedron. **B.** Structure – first flat position. **C.** Structure - second flat position. **D.** $P_0$ Cap Construction Parameters. **E.** Extension Parameters.

In this hecatohedron construction (Fig. 14), the intermediate octahedra are created by scaling the octahedron $P_0$ with factors that control the length of the extension of the length of edges $|A_iA_{i+1}|$. The octahedron $P_0$ is created using the Eqs. 2a and 3a with the constant k as defined in Eq. 2.

Table VI contains construction parameters for the other polyhedra used in the body of this paper.

| Figs. 5 and 6. | $|X_0A_0|$=10, $|X_0B_0|$=10.5, $|X_0C_0|$=9.5, $|X_0D_0|$=9, $|A_0B_0|$=1.5, $|B_0C_0|$=1.75. $|C_0C_1|$=5. |
|---|---|
| Fig. 7A | $f_2$=2.5, $f_3$=4, $f_4$=5, $f_5$=6. |
| Fig. 8 | $|X_0A_0|$=10, $\angle A_0X_0B_0$=17°, $\angle B_0X_0C_0$=47°, $\angle B_0A_0X_0$=65°, $\angle C_0B_0X_0$=40°, $|B_0B_1|$=0, $|C_0C_1|$=0. $|D_0D_1|$=0. |
| Fig. 9 | $|X_0A_0|$=10, $\angle A_0X_0B_0$=17°, $\angle B_0X_0C_0$=47°, $\angle B_0A_0X_0$=65°, $\angle C_0B_0X_0$=40°, $|B_0B_1|$=1.5, $|C_0C_1|$=0. $|D_0D_1|$=0. |

**Table VI.** Example Polyhedra Parameters.

The hexadecahedron that is used to illustrate the construction of tori (Figs. 5 and 6) is a sub-type I-OEE construction parameterized as shown in Table VI. These parameters have been selected to create a relatively elongated cap in order to better illustrate subsequent torus construction.

The structure (Fig. 7A) is based upon the $P_0$ octahedron for the hexadecahedron (Fig. 5) with scale factors as shown in Table VI applied to define five distinct polyhedra that appear in the structure. Instances of the $P_0$ octahedron, which has the implied scale factor $f_1$=1, appear along the bottom of Fig, 7A.